\documentclass[12pt]{amsart}
\usepackage[all,cmtip]{xy}

\usepackage{amssymb}
\usepackage{latexsym}
\newtheorem{theorem}{Theorem}[section]
\newtheorem{lemma}[theorem]{Lemma}
\newtheorem{corollary}[theorem]{Corollary}
\newtheorem{proposition}[theorem]{Proposition}

\newtheorem{definition}[theorem]{Definition}

\newcommand{\ble}{\begin{lemma}}
\newcommand{\ele}{\end{lemma}}
\newcommand{\bth}{\begin{theorem}}
\renewcommand{\eth}{\end{theorem}}
\newcommand{\bpr}{\begin{proposition}}
\newcommand{\epr}{\end{proposition}}
\newcommand{\bco}{\begin{corollary}}
\newcommand{\eco}{\end{corollary}}
\newcommand{\bde}{\begin{definition}}
\newcommand{\ede}{\end{definition}}
\newcommand{\beq}{\begin{equation}}
\newcommand{\eeq}{\end{equation}}
\newcommand{\bpf}{\begin{proof}}
\newcommand{\epf}{\end{proof}}

\newcommand{\ben}{\begin{enumerate}}
\newcommand{\een}{\end{enumerate}}
\newcommand{\bcon}{\begin{conj}}
\newcommand{\econ}{\end{conj}}
\newcommand{\bex}{\begin{exa}}
\newcommand{\eex}{\end{exa}}
\newcommand{\barr}{\begin{array}}
\newcommand{\earr}{\end{array}}
\newcommand{\btab}{\begin{tabular}}
\newcommand{\etab}{\end{tabular}}
\newcommand{\bea}{\begin{eqnarray*}}
\newcommand{\eea}{\end{eqnarray*}}
\newcommand{\bce}{\begin{center}}
\newcommand{\ece}{\end{center}}
\newcommand{\bpi}{\begin{picture}}
\newcommand{\epi}{\end{picture}}
\newcommand{\bfi}{\begin{figure} \begin{center}}
\newcommand{\efi}{\end{center} \end{figure}}

\newcommand{\bsl}{\begin{slide}{}}
\newcommand{\esl}{\end{slide}}

\newcommand{\hso}[1]{\hspace{-1pt}}

\newcommand{\sbe}{\subseteq}

\newcommand\nid{{\noindent}}
\newcommand\cp{{\mathcal{P}}}
\newcommand\cl{{\mathcal{L}}}
\newcommand\cx{{\mathcal{X}}}

\def\<{\langle}
\def\>{\rangle}

\newcommand{\De}{\Delta}

\newcommand{\cA}{{\mathcal A}}

\newcommand{\cD}{{\mathcal D}}

\newcommand{\cE}{{\mathcal E}}

\newcommand{\cH}{{\mathcal H}}

\newcommand{\cL}{{\mathcal L}}

\newcommand{\cP}{{\mathcal P}}

\newcommand{\cS}{{\mathcal S}}

\newcommand{\cW}{{\mathcal W}}

\newcommand{\xx}{\mathsf{x}}
\newcommand{\yy}{\mathsf{y}}

\newcommand{\hcS}{{\widehat{\mathcal S}}}

\renewcommand{\bar}{\overline}

\newcommand{\id}{\mathop{\rm id}\nolimits}

\newcommand{\Rad}{\mathop{\rm Rad}\nolimits}

\newcommand{\cho}{\choose}

\def\flexbox#1{\mathchoice{\mbox{#1}}{\mbox{#1}}{\mbox{\scriptsize #1}}%
{\mbox{\tiny #1}}}
\def\SL{\mathop{\flexbox{\rm SL}}}

\newcommand{\SU}{\mathop{\flexbox{\rm SU}}}
\newcommand{\SG}{\mathop{\flexbox{\rm SG}}}

\newcommand{\PG}{\mathop{\flexbox{\rm PG}}}

\def\pset{{\mathcal P}}
\def\lset{{\mathcal L}}
\def\cD{{\mathcal D}}
\newcommand{\Char}{\mathop{\flexbox{\rm Char}}}

\newcommand{\map}{\longrightarrow}

\newcommand{\FF}{{\mathbb F}}
\newcommand{\KK}{{\mathbb K}}
\newcommand{\NN}{{\mathbb N}}
\newcommand{\PP}{{\mathbb P}}

\newcommand{\sform}{{\sf s}}
\newcommand{\hform}{{\sf h}}
\newcommand{\form}{{\sf f}}

\newcommand{\gr}{{\rm gr}}%

\newcommand{\GR}{\mbox{Gr}}

\newcommand{\dfn}{\em}
\newcommand{\Res}{\mathop{\rm Res}\nolimits}

\newcommand{\proj}{\mathop{\rm proj}}%
\newcommand{\dist}{{\rm d}}%

\newcommand{\Sp}{\mathop{\rm Sp}}

\newcommand{\wA}{{{}^2\hspace{-.2em}A}}

\newcommand{\AI}[1]{\item[\rm{(#1)}]}

\newcommand{\emb}{\mathsf{e}}

\newcommand{\what}{\widehat}
\newtheorem{thm}{Theorem}[section]

\newtheorem{conj}[thm]{Conjecture}
\newtheorem{exa}[thm]{Example}

\newenvironment{frontmatter}{\pagestyle{empty}}{\pagestyle{plain}}
\newenvironment{keyword}{\noindent Keywords:}{}
\newcommand{\ead}{\email}

\newcommand{\corauthref}[1]{}
\newcommand{\sep}{\ }

\begin{document}
\begin{frontmatter}



\title{The generating rank of the  unitary and symplectic Grassmannians}


\author{Rieuwert J. Blok\corauthref{cor}}
\address{Department of Mathematics and Statistics,
Bowling Green State University,
Bowling Green, OH 43403,
U.S.A.}
\ead{blokr@member.ams.org}

\author{Bruce N. Cooperstein}
\address{Department of Mathematics,
University of California,
Santa Cruz,
U.S.A.}
\ead{coop@ucsc.edu}

   \begin{abstract}
We prove that the Grassmannian of totally isotropic $k$-spaces of the
polar space associated to
 the unitary group $\SU_{2n}(\FF)$  ($n\in\NN$) has generating rank ${2n
\choose k}$ when $\FF\ne \FF_4$.
We also reprove the main result of Blok~\cite{Bl2007}, namely that the
Grassmannian of totally isotropic $k$-spaces
 associated to the symplectic group $\Sp_{2n}(\FF)$ has generating rank
${2n\choose k}-{2n\choose k-2}$, when $\Char(\FF)\ne 2$.
    \end{abstract}
\maketitle
\begin{keyword}
unitary group\sep Hermitian form\sep symplectic group\sep Grassmannian\sep generating
rank\sep embedding\\
AMS subject classification (2000):
    Primary 51A50;
    Secondary 51A45, 05A15.
\end{keyword}
\end{frontmatter}

\section{Introduction}\label{sec:introduction}
Generating sets of point-line geometries serve both theoretical and
computational purposes.
For instance, if the geometry admits a finite-dimensional absolutely
universal embedding and if the generating rank equals the (vector)
dimension of a particular embedding, then the latter is absolutely
universal.
On the other hand, minimal generating sets may serve in creating
computer models of point-line geometries.

Until now, the projective building of type $A_n$ associated to the group $\SL_{n+1}(\FF)$ and the symplectic
building of type $C_n$ associated to $\Sp_{2n}(\FF)$ ($\Char(\FF)\ne 2$)
are the only ones of which the generating rank is known for each of its
$k$-shadow spaces ($k$ a single node of the diagram). (see Cooperstein and Shult~\cite{CoSh1997}, Blok and
Brouwer~\cite{BlBr1998a}, Blok~\cite{Bl2003,Bl2007},
Cooperstein~\cite{Co1998,Co2003}, and De Bruyn and Pasini~\cite{DePa2007}.)

In this paper we deal with the unitary building of type
$\wA_{2n-1}(\FF)$ associated to the group
 $\SU_{2n}(\FF)$ over a field $\FF\not\cong\FF_4$ and exhibit minimal
generating sets for all associated unitary $k$-Grassmannians. For a few
of these geometries a minimal generating set has already been found.
Since the polar space associated to $\SU_4(q^2)$ has more points per line than lines per point, it is generated by the four points of an apartment.  A result by Blok and Brouwer~\cite[Theorem 2.1]{BlBr1998a} then implies that the generating rank of the unitary polar space associated to the group
$\SU_{2n}(\FF_{q^2})$ is $2n$.

In~\cite{Co1997} Cooperstein proved that the dual polar space associated
to $\SU_{2n}(q^2)$ has generating rank
 ${2n \choose n}$ when $q^2>2$.
 This result was generalized to include infinite fields by De Bruyn and
Pasini~\cite{DePa2007}.
It was proved by A.~E.~Brouwer that the generating rank is 
 at least $\frac{4^n+2}{3}$ when $q^2=4$. His conjecture that equality holds was confirmed by P. Li~\cite{Li2002}.

We shall employ these results to obtain minimal generating sets for all
$k$-shadow spaces of this building in a unified way. 
\bth\label{genranksugrassmannians}
 The $k$-Grassmannian of the polar space associated to
 $\SU_{2n}(\FF)$ ($\FF$ a field) has generating rank ${2n \choose k}$ if
$\FF \ne \FF_4$ or $k=1$.
\eth

The techniques we
use also allow us to describe
 generating sets for all symplectic $k$-Grassmannians over fields of
characteristic different from $2$ in almost exactly the same manner,
thus giving an alternative proof for the main result in Blok~\cite{Bl2007}.

\bth\label{genrankspgrassmannians}{\rm (Blok~\cite{Bl2007})}
 The $k$-Grassmannian of the polar space associated to
 $\Sp_{2n}(\FF)$ ($\FF$ a field) has generating rank ${2n \choose k}-{2n
\choose k-2}$ if $\Char(\FF)\ne 2$.
\eth

More on generating sets can be found in Cooperstein~\cite{Co2003} and
Blok~\cite{Bl1999}.

In Section~\ref{section:preliminaries} we define some basic concepts, including the generating rank of a point-line geometry.
We also introduce all geometries under study. We then prove 
Proposition~\ref{prop:Gamma_1 spanned by apartment}, which is a special case of Theorems~\ref{genranksugrassmannians}~and~\ref{genrankspgrassmannians}, and present Theorem~\ref{thm:bruce's theorem}, which is a special case of Theorem~\ref{genranksugrassmannians}. These two will form the basis of an inductive proof of our main results.
In Section~\ref{section:emb} we describe embeddings of $\Gamma_k^\form$
of dimension $d_k$, where
 $$d_k=\left\{\begin{array}{@{}ll}
 {2n\choose k} & \mbox{if }$\form=\hform$\\
 {2n\choose k}-{2n\choose k-2}&\mbox{if } \form=\sform\\
 \end{array}\right.
 $$
This section also contains a short proof of the probably well-known result that the $\SU(V)$ module $\bigwedge^k V$ is irreducible when $k<\dim(V)$ (see Theorem~\ref{thm:exterior power is SU-irreducible}).

It then suffices (see~(\ref{eq:genrank})) to show that
$\Gamma_k^\form$ has a generating set of size exactly
 $d_k$.
We do this in Section~\ref{section:gensets}.
The definition of the generating set, as well as the proof that the
proposed set generates the geometry and has the desired
 cardinality is inductive.
Two inductive tools are described in Section~\ref{section:ext}.
%
\section{Preliminaries}\label{section:preliminaries}
A {\em point-line geometry} is a pair $\Gamma=(\pset,\lset)$ where
$\pset$ is a
set whose elements are called `points' and $\lset$ is a collection of
subsets of
$\pset$ called `lines' with the property that any two distinct points
belong to at most
one line.
If $\pset$ and $\lset$ are not mentioned explicitly, the sets of points
and lines
 of a point-line geometry $\Gamma$ are denoted $\pset(\Gamma)$ and
$\lset(\Gamma)$.

A {\em subspace} of $\Gamma$ is a subset $X\subseteq \pset$ such that
any line
containing at least two points of $X$ is itself entirely contained in $X$. We call $X$ {\em proper} if $X\subsetneq \pset$.
A {\em hyperplane} of $\Gamma$ is a proper subspace that meets every line.

\paragraph{Projective embeddings and generating sets}
\noindent
The {\em span} of a set $\cS\subseteq \pset$ is the smallest subspace of $\Gamma$
containing
$\cS$; it is the intersection of all subspaces containing $\cS$ and is
denoted by
$\langle \cS\rangle_\Gamma$.
We say that $\cS$ is a {\em generating set} (or {\em spanning set}) for
$\Gamma$
 if $\langle \cS\rangle_\Gamma=\pset$.

For a vector space $W$ over some field $\FF$, the {\em projective geometry}
 associated to $W$ is the point-line geometry  $\PP(W)=(\pset(W),\lset(W))$
 whose points are the $1$-spaces of $W$ and whose lines are the sets of
$1$-spaces
 contained in some $2$-space.

A {\em projective embedding} of a point-line geometry
 $\Gamma=(\pset,\lset)$ is a pair $(\epsilon, W)$, where $\epsilon$
 is an injective map $\pset\stackrel{\epsilon}{\map}\pset(W)$
 that sends every line of $\lset$ onto a line of $\lset(W)$, and with
the property
 that $$\langle\epsilon(\pset)\rangle_{\PP(W)}=\pset(W).$$
In the literature, this is sometimes referred to as a {\em strong} or {\em full} projective
embedding.
The {\em dimension} of the embedding is the dimension of the vector
space $W$.
In this paper we will assume both $\dim(W)$ and $|\cS|$, ($\cS$ as
above) to be finite. Then, since $\epsilon(\langle
\cS\rangle_\Gamma)\sbe \langle \epsilon(\cS)\rangle_W$, for any
generating set $\cS$ and any embedding
 $(\epsilon,W)$ we have
\beq\label{eq:genrank}
\dim(W)\le |\cS|.
\eeq
In case of equality, $\cS$ is a minimal generating set and $(\epsilon,W)$ is an
embedding of
 maximal dimension.
We then call $\dim(W)=|\cS|$ the {\em generating rank} of $\Gamma$.
\medskip

We briefly describe the particular geometries we will discuss in this paper.
\paragraph{The Projective Grassmannian}
\noindent
Let $V$ be a vector space over some field $\FF$.
For any $k$ with $1\le k\le \dim(V)-1$, the {\em projective
$k$-Grassmannian} associated to $V$
 is the point-line geometry $\GR(V,k)$ whose points are the $k$-spaces
of $V$
 and whose lines are sets of the form
 $$\{K\mbox{ a $k$-space in $V$}\mid  D\subseteq K\subseteq U\}$$
 for some $(k-1)$-space $D$ and $(k+1)$-space $U\supseteq D$.
One verifies that the map
 \beq\label{eqn:projective grassmann embedding}
 \begin{array}{rl}
 \emb_k\colon \GR(V,k)&\to \PG(\bigwedge^k V)\\
 U & \mapsto \bigwedge^k U
 \end{array}
 \eeq
defines a full projective embedding for $\GR(V,k)$, called the {\em Grassmann embedding}. 
Grassmann embeddings have been studied in connection with 
many topics, including hyperplanes~\cite{Sh1992}, highest weight modules~\cite{Blo2010}, algebraic varieties~\cite{CohCoo1998} and Schubert Calculus~\cite{Fu1997}.

\medskip

\paragraph{The Unitary and Symplectic Grassmannians}
\noindent
Let $V$ be a vector space of dimension $2n$ over the field $\FF$ endowed
with a form $\form$.
We assume that either $\form$ is non-degenerate and symplectic, or
 $\form$ is non-degenerate $\sigma$-Hermitian of Witt index $n$.
This means that $\sigma$ is the generator of the Galois group
of a quadratic field extension $\FF/\mathbb{K}$, which is assumed to be separable in the case the characteristic is two.
The norm $N_\sigma\colon \FF\to\KK$ is then trace-valued.
We sometimes write $\KK=\FF^\sigma$ to emphasize that $\KK$ is the fixed subfield of $\FF$ under $\sigma$.
In the respective cases we sometimes write $\sform$ or $\hform$ for $\form$.
For a subspace $U$ of $V$ we define
$$U^\perp=\{v\in V\mid \form(u,v)=0\, \forall u\in U\}.$$
We write $U\perp W$ if $W\sbe U^\perp$.
The {\dfn radical} of a subspace $W$ is
$$\Rad(W,\form)=W^\perp\cap W.$$
A subspace $U$ of $V$ is called {\em totally isotropic} (t.i.) with
respect to the
 form $\form(\cdot,\cdot)$ if $U\sbe U^\perp$.
It is called non-degenerate if $\Rad(U,\form)=\{0\}$.

The {\dfn polar building} $\Gamma$ associated to $\form$ is the
spherical building whose objects of type $i\in I=\{1,2,\ldots,n\}$ are
the $i$-spaces in $V$ that are totally isotropic  with respect to $\form$.
Two objects are incident when one contains the other as a subspace of $V$ (see e.g.~\cite{BuSh1974,Ti1974,Vel1959}).

The {\dfn polar space} associated to $\form$ is the point-line geometry
$\Gamma^\form_1$ whose
 points are the t.i.\ $1$-spaces of $V$ and whose lines are sets of
$1$-spaces of the form
  $$\{P \mbox{ a t.i.\ $1$-space of $V$}\mid P\sbe L\}$$ for some t.i.\
$2$-space $L$.
We sometimes call t.i.\ $3$-spaces {\em planes}.
A {\dfn hyperbolic line} $H$ is a $2$-space in $V$ such that $\form$
restricted to $H$ is non-degenerate of Witt index $1$.
Since both the symplectic and unitary polar space under study in this paper have Witt index $n=\frac{1}{2}\dim(V)$, the vector space $V$ has a basis $\cE=\{e_i,e_{n+i}\mid i=1,2,\ldots,n\}$ that is {\em hyperbolic} with respect to $\form$. That is, 
 $\form(e_i,e_j)=0$ for all $1\le i\le j\le 2n$ with $j\ne n+i$
  and $\form(e_i,e_{n+i})=1$ for all $i=1,2,\ldots,n$.
Each hyperbolic basis $\cE$ for $V$ gives rise to an apartment 
 $\Sigma(\cE)$ of the polar building; it can conveniently be described as the collection $\Sigma(\cE)=\{E_{J,K}\mid J,K\}$, where 
 $E_{J,K}=\langle e_j,e_{k+n}\mid j\in J, k\in K\}$ and $(J,K)$ runs over all pairs of subsets of $\{1,2,\ldots,n\}$ such that
 $J\cap K= \emptyset\ne J\cup K$.
 
The {\em (polar) $k$-Grassmannian} associated to $\form$ is the
point-line geometry $\Gamma^\form_k(V)=(\cP_k,\cL_k)$
 whose points are the t.i.\ $k$-spaces and whose lines are the sets of
the form
 $$\{K\mbox{ a t.i.\ $k$-space in $V$}\mid  D\subseteq K\subseteq U\}$$
 for some t.i.\ $(k-1)$-space $D$ and t.i.\ $(k+1)$-space $U\supseteq D$.
In case $k=n$, the lines are of the form
 $$\{K\mbox{ a t.i.\ $k$-space in $V$}\mid  D\subseteq K\}$$
for some t.i.\ $(n-1)$-space $D$.
Whenever $V$ or $\form$ is clear from the context, we'll drop it from
the notation. Note that $\Gamma^\form_k(V)$ is a subgeometry of $ \GR(V,k)$, i.e.\ $\cP(\Gamma^\form_k)\sbe \cP(\GR(V,k))$ and
 $\cL(\Gamma^\form_k)\sbe \cL(\GR(V,k))$.

The {\dfn $k$-residue} or simply {\dfn residue in $\Gamma_k$} of an
object $X$ of $\Gamma$,
 is denoted $\Res_k(X)$ and is the subgeometry of points and lines of
$\Gamma_k$ that are incident to $X$
  when viewed in $\Gamma$.

We denote the group of linear transformations of determinant $1$
preserving the form $\form$ on $V$ (i.e.\ the symmetry group of $\form$) by
 $\SG(V,\form)$.
Thus, $\SG(V,\sform)=\Sp(V)\cong \Sp_{2n}(\FF)$ and
 $\SG(V,\hform)=\SU(V)\cong \SU_{2n}(\FF)$.
Clearly $\SG(V,\form)$ is an automorphism group of $\Gamma_k^\form$ for
all $k=1,2,\ldots,n$ in that
 it preserves points and lines and the incidence between them.

The following special case of Theorems~\ref{genranksugrassmannians}~and~\ref{genrankspgrassmannians} will be instrumental in proving those theorems. 
\bpr\label{prop:Gamma_1 spanned by apartment}
Let $\Gamma_1$ be the polar space associated to a non-degenerate
symplectic or Hermitian form $\form$ of Witt index $n$ on a vector space
$V$ of dimension $2n$ over $\FF$. Moreover, in case the form is
symplectic assume that $\Char(\FF)\ne 2$.
Then $\Gamma_1$ has generating rank $2n$ and it is generated by the $2n$
points of an apartment i.e.\ a hyperbolic basis.
\epr

\bpf
Theorem 2.1 of~\cite{BlBr1998a} shows that a non-degenerate polar space of rank $n\ge 2$ is generated by the $2n$ points of an apartment if this is so in rank $2$.
For the symplectic rank-$2$ case this is proved in e.g.\ loc.~cit.~or~\cite{CoSh1997}.
For the Hermitian rank-$2$ case this is the content of Lemma~\ref{lem:SU4} below.
\epf

\ble\label{lem:SU4}
The unitary rank $2$ polar space associated to $\SU_4(\FF)$ is spanned
by four points.
\ele
\bpf
\nid  Let  $\FF$, $\KK$ and $\sigma$ be as defined above.
For convenience denote images under
$\sigma$ using the bar notation. Assume $\FF$ is generated over
$\mathbb{K}$ by the element $\delta$ which we further assume satisfies
$\delta + \overline{\delta} = 0$ if the characteristic is not two.

\medskip

\noindent 
We assume that $\Gamma_1=(\cP,\cL)$ is the polar space associated to the  $\sigma$ Hermitian form $\hform\colon V\times V\to \FF$, where $V$ is a 4-dimensional space over $\FF$.
Let $\cE=\{e_i,e_{2+i}\mid i=1,2\}$ be a hyperbolic basis for $V$ with respect to $\hform$.
%
%
%
We
will  prove that the set $\cS = \{\langle e_1\rangle, \langle e_2
\rangle, \langle e_3 \rangle, \langle e_4\rangle\}$ is a generating
set for $(\cp, \cl).$  Before proving this we  introduce some notation.

\medskip

\nid 
Let  $L = \langle e_1, e_2\rangle$, $M = \langle
e_3, e_4\rangle$,  $P_a = \langle e_1 + ae_2\rangle$, and $Q_a=\langle -\overline{a}e_3 + e_4\rangle.$  

\medskip

\nid 
The following well-known facts about $\Gamma_1$ and  the isometry group $\SU_4(\FF)$ are easy to prove:
\renewcommand{\theenumi}{\roman{enumi}}
\begin{enumerate}
\item\label{isotropic} If $\xx, \yy$ are isotropic vectors and $\hform(\xx, \yy) = 1$ then  a
vector $\xx + a\yy$ is isotropic if and only if $trace(a) = a +
\overline{a} = 0.$ When the characteristic is two this implies that $a
\in \mathbb{K}$ and otherwise $a \in \delta \mathbb{K}.$  


\item\label{PaQa} The unique point on $M$ which is orthogonal to $P_a$ is $Q_a.$


\item\label{tau} If $\tau \in \SL(V)$ and the matrix of $\tau$ with respect to
${\mathcal{B}}$ is
$$\left(\begin{array}{@{}cc@{}} 
A & 0_{22}\\
0_{22} & B\end{array}\right),$$
where $A, B$ are $2 \times 2$ matrices,  then $\tau\in \SU_4(\FF)$ if and only if $\overline{B} = A^{-T}.$


\item\label{trans} A consequence of (\ref{tau}) is that $\SU_4(\FF)$ is
transitive on the points of $L.$


\item\label{nonorth trans} A further consequence of (\ref{tau}) is that $\SU_4(\FF)$ is transitive on pairs $(P,Q)$ where $P \in L$, $Q\in M$ and $P$ and $Q$
are not orthogonal.
\end{enumerate}

\nid We now prove our assertion about the generation of $\cp.$    Let
$\cx$ denote the subspace of $\cp$ generated by $\cS.$  Let $R$ be an
arbitrary point of $\cp.$  If $R \in L \cup M$ there is nothing to
prove.   Otherwise, there are points $P \in L, Q \in M$ such that $R
\in \langle P, Q\rangle.$  If $P \perp Q$ then $R$ is in $\cx$  since
both $P$ and $Q$ are.  Thus, we may assume that $P$ and $Q$ are not
orthogonal.  By (\ref{nonorth trans}) we can assume that $P = \langle e_1\rangle$ and $Q =
\langle e_3\rangle.$   By (\ref{isotropic}) we may therefore assume that $R = \langle
e_1 + c\epsilon e_3 \rangle$ for some $c \in \mathbb{K}$ where $\epsilon =
1$ if the characteristic is two and $\epsilon = \delta$ otherwise.

\medskip

\nid As previously pointed out, the points $P_a, Q_a$ are orthogonal and
in $\cx$ and therefore the line $\langle P_a, Q_a\rangle$ is contained
in $\cx.$  In particular, the point 
$D_a  = \langle (e_1 + ae_2)
-(c\epsilon/ \overline{a})(-\overline{a}
e_3 + e_4)\rangle =
\langle e_1 + ae_2 + c\epsilon e_3 -
(c\epsilon/\overline{a})
e_4\rangle$ is isotropic and in $\cx.$
The same holds for 
$D_b = \langle e_1 + be_2 + c\epsilon e_3 -
(c\epsilon/\overline{b})
e_4\rangle$.

\nid 
Using that $\bar{\epsilon}=-\epsilon$ and $\bar{c}=c$, we find
$$\hform(e_1 + ae_2 + c\epsilon e_3 -
\frac{c\epsilon}{\overline{a}}e_4, e_1 + be_2 + c\epsilon e_3 -
\frac{c\epsilon}{\overline{b}}e_4) =$$
$$ -c\epsilon + \frac{c\epsilon a}{b} + c\epsilon -\frac{c\epsilon
\overline{b}}{\overline{a}} = c\epsilon(\frac{a}{b} -
\frac{\overline{b}}{\overline{a}})$$ 
so that $D_a\perp D_b$ if and only if $a\overline{a} =
b\overline{b}.$

\nid
It is possible to choose distinct $a$ and $b$  such that 
$a\overline{a} = b\overline{b}.$ 
We claim that $R=\langle e_1 + c\epsilon e_3\rangle \in \langle D_a, D_b\rangle$. Namely, we have 
$$b(e_1 + ae_2 + c\epsilon e_3 -
\frac{c\epsilon}{\overline{a}}e_4) - a(e_1 + be_2 + c\epsilon e_3
- \frac{c\epsilon}{\overline{b}}e_4) = $$
$$(b - a)e_1 + (b - a)c\epsilon e_3 +
(c\epsilon)\left[\frac{a}{\overline{b}} - \frac{b}{\overline{a}}\right] e_4=$$
$$ (b-a)[e_1 + c\epsilon e_3].$$

\medskip

\nid We can conclude from $R \in \langle D_a, D_b\rangle, D_a \perp
D_b$, and $D_a, D_b \in \cx$ that $R \in \cx$  as required. 
\epf

We shall also need the following result from De Bruyn and
Pasini~\cite{DePa2007} and Cooperstein~\cite{Co1998,Co1997} on the
generating rank of $\Gamma_n$.
\bth\label{thm:bruce's theorem}
Let $\Gamma_n$ be the dual polar space associated to a non-degenerate
symplectic or Hermitian form $\form$ on a vector space $V$ of dimension
$2n$ over $\FF$.
Then
\begin{itemize}
\AI{a} $\Gamma_n$ has generating rank ${2n\cho n}$ if $\form=\hform$ is
Hermitian and $\FF\ne \FF_4$, and
\AI{b} $\Gamma_n$ has generating rank ${2n\cho n}-{2n\cho n-2}$ if
$\form=\sform$ is symplectic and $\FF\not\cong\FF_2$.
\end{itemize}
\eth

\section{Embeddings}\label{section:emb}
We describe the fairly well-known Grassmann embeddings of the polar
Grassmannian $\Gamma_k^\form$.
We consider arbitrary $(\form,k)$ with $1\le k\le n$ and
$\form=\sform,\hform$, excluding - for the moment - only the case
$(\hform,n)$.
We then define the following map on the point-set of $\Gamma_k^\form$:
\beq\label{eqn:embedding map}
\hspace{1.6in}
\begin{array}{rl}
\emb_\gr\colon\Gamma_k& \to \PG(\bigwedge^k V)\\
U & \mapsto \bigwedge^k U\\
\end{array}
\eeq
Here for any ordered set of vectors $\{u_1,\ldots,u_k\}$ define
$\bigwedge^k \{u_1,\ldots,u_k\}= u_1\wedge u_2\wedge \cdots\wedge u_k$.
One verifies that, if $\{u_1,\ldots,u_k\}$ is a basis for a
subspace $U$, then the $1$-space $\bigwedge^kU=\langle
\bigwedge^k\{u_1,\ldots,u_k\}\rangle$ only depends on
  $U$, not on the chosen ordered basis.
Thus, $\emb_\gr$ is simply the restriction to $\Gamma_k$ of the Grassmann embedding defined in~(\ref{eqn:projective grassmann embedding}). Since the lines of $\Gamma_k$ are lines of 
 $\GR(V,k)$, it follows that $\emb_\gr$ is in fact an embedding of
$\Gamma_k$ into $\PG(V_\gr)$ for some subspace $V_\gr$ of $\bigwedge^k V$.

\bth\label{thm:exterior power is SU-irreducible}
Let $V$ be a vector space of dimension $n$ over $\FF$ endowed with a non-degenerate Hermitian form.
Then, the vector space $W=\bigwedge^k V$ is  irreducible as a module for the unitary group $\SU_n(\FF)=\SU(V)$.
\eth
\bpf
We shall prove that the space $W=\bigwedge^k V$, which has dimension
${2n\choose k}$ over $\FF$ is an irreducible module for the unitary group
 $\SU_n(\FF)=\SU(V)$.
We first note that $\SU(V)\le\SL(V)$ and that $W$ is irreducible for $\SL(V)$ over any field.
This is well-known, but it is also easy to see.
View $\SL(V)\cong\SL_n(\FF)$ via a basis $\{a_1,\ldots,a_n\}$ for $V$ and let $H$ be the group of diagonal martrices.
Then since $H$ is abelian, any submodule $U$ of $W$ is the direct sum of weight spaces, i.e.\ common eigenspaces for all elements of $H$.
One verifies that since $|\FF|\ge 4$ the decomposition of $W$ into weight spaces for $H$ is unique and equals $\bigoplus_{w\in \cW} w$, where $\cW=\{\langle a_{i_1}\wedge \cdots\wedge a_{i_k}\rangle_W\mid 1\le i_1<\cdots<i_k\le n\}$. Thus, if $U$ is non-zero it contains at least one $1$-space from $\cW$, but since the subgroup $N\le \SL_n(\FF)$ of monomial matrices is transitive on $\cW$ and $\cW$ spans $W$, we find  $U=W$.
 
In view of what we just saw,  it suffices to show that if a subspace of $W$ is $\SU(V)$-invariant, then it is $\SL(V)$-invariant.

To this end we examine the inclusion $\SU(V)\le \SL(V)$ a little closer.
Each root group $T^\sigma$ of the former is contained in
a unique root group $T$ of the latter.
Namely, there is an $\hform$-isotropic vector $u$, such that
$$\begin{array}{lll}
T^\sigma&=T_u^\sigma&=\langle t_{u,\lambda}\mid
\lambda\in\FF^\sigma\rangle,\\
T&=T_u&=\langle t_{u,\lambda}\mid \lambda\in\FF\rangle,\\
\end{array}$$
where $t_{u,\lambda}(v)=v+\lambda\delta\hform(v,u)u$ and 
 $\FF$ is generated over $\FF^\sigma$ by $\delta$, and moreover, 
  $\delta+\delta^\sigma=0$ whenever $\Char(\FF)\ne 2$. 
One can show that $\SU(V)$ is generated by all such groups $T_u^\sigma$
 (See e.g.\ \cite{Di1955,Ta1992}).
What is more, we claim that:
\beq\label{eqn:SL gen T_u}
\SL(V)=\langle T_u \mid u\in V\mbox{ $\hform$-isotropic}\rangle.
\eeq
Proof of (\ref{eqn:SL gen T_u}):
For $n=2$ it is easy to verify that any two non-orthogonal isotropic
vectors $e$ and $f$ suffice.
For $n\ge 3$, let $a_1,\ldots,a_n$ be a basis for $V$ that is orthogonal
with respect to $\hform$. Then
\begin{eqnarray}
\SL(V)&=\langle \SL(\langle a_i,a_{i+1}\rangle)\mid 1\le i\le n-1\rangle.&\label{eqn:SLV SL2} 
\end{eqnarray}
Presenting $\SL(V)\cong \SL_n(\FF)$ via the basis $\{a_1,\ldots,a_n\}$, it is rather easy to show, using commutators of elementary matrices, that the right hand-side of  (\ref{eqn:SLV SL2}) contains all elementary matrices, so that equality follows.

Since the norm $N_\sigma\colon \FF\to \FF^\sigma$ is surjective, for each $i$ there exist two isotropic vectors
$e_i$ and $f_i$ such that $\langle a_i,a_{i+1}\rangle=\langle
e_i,f_i\rangle$ (see e.g.~\cite[Ch. 10]{Ta1992}).
Now (\ref{eqn:SL gen T_u}) follows by applying the $n=2$ case
$n-1$ times.

Next, we compute the action of $(t_{u,\lambda}-\id)$ on $W$.
Select a basis $u=v_1,\ldots,v_n$ for $V$ such that $\hform(v_i,v_1)=0$ for all $i\ne 2$ and $\delta\hform(v_2,v_1)=1$.
The set of all pure vectors $\wedge_{j=1}^k v_{i_j}$ with $1\le i_1<i_2<\cdots<i_k\le n$
 forms a basis for $W$.
For such a basis we have 
 $$\begin{array}{rl}
 (t_{u,\lambda}-\id)(\wedge_{j=1}^k v_{i_j})&=
   \wedge_{j=1}^k t_{v_1,\lambda}(v_{i_j})-\wedge_{j=1}^k v_{i_j} \\
  &=\sum_{j=1}^k  v_{i_1}\wedge\cdots\wedge v_{i_{j-1}}\wedge \lambda\delta\hform(v_{i_j},v_1)v_1\wedge v_{i_{j+1}}\wedge\cdots\wedge v_{i_k}\\
  \end{array}$$
since any pure wedge product with two or more occurrences of $v_1$ is the zero vector.
Using that $\hform(v_i,v_1) = 0$ for all $i\ne 2$, we find that
$$ (t_{u,\lambda}-\id)(\wedge_{j=1}^k v_{i_j})=
 \left\{\begin{array}{ll}
 \lambda v_1\wedge v_{i_2}\wedge\cdots\wedge v_{i_k} & \mbox{ if }i_1=2,\\
 0 & \mbox{ else. }\\
 \end{array}\right..$$
Since the image of any basis vector is either $0$ or again a basis vector multiplied by $\lambda$, it follows that for any $1$-dimensional subspace $A$ of 
 $W$ we have
  $[A,T_u]=[A,T_u^\sigma]$.
 (Recall that if a group $T'$ acts on a vector space $V'$, then for any subspace $A'\le V'$ one writes 
  $[A', T']=\langle (t'-\id)(a)\mid t'\in T'\ a\in A'\rangle_{V'}$. See e.g.~\cite{KuSt2004}.) 
Clearly an arbitrary subspace $U$ of $W$ is invariant under $T_u^\sigma$ if and only if $[A,T_u^\sigma]\le U$ for any $1$-dimensional subspace $A$ of $U$, and the same holds for $T_u$.
Thus since $[A,T_u^\sigma]=[A,T_u]$, for all such $A$,  $U$ is $T_u^\sigma$-invariant if and only if it is $T_u$-invariant.
Our claim now follows from (\ref{eqn:SL gen T_u}).
\epf

\bco\label{cor:hermitian grassmann embedding}
Let $1\le k\le n-1$.
Then, the Grassmann embedding $(\emb_\gr,V_\gr)$ for $\Gamma_k^\hform$
has dimension ${2n\choose k}$.
More precisely, we have $V_\gr=\bigwedge^k V$.
\eco
\bpf
As $\SU(V)$ acts transitively on the t.i.\ $k$-spaces in $V$, the  space $V_\gr$ is an $\SU(V)$-submodule of $\bigwedge^k V$.  By Theorem~\ref{thm:exterior power is SU-irreducible}, $V_\gr=\bigwedge^k V$.
\epf

\medskip
\noindent
The Grassmann embedding for $\Gamma_n^\hform$ is the same map $e_\gr$ as
given in (\ref{eqn:embedding map}), but it has the property
that the image $(\bigwedge^nV)^\sigma$ can be viewed as a Baer subspace
over the fixed field $\FF^\sigma$.
Again, $(\bigwedge^nV)^\sigma$, as a vector space over $\FF^\sigma$ has
dimension ${ 2n\choose n}$, but its vectors do span $\bigwedge^n V$ (see
e.g.~\cite{De-2008a}).
In this paper we will call $\emb_\gr=\emb_n^\sigma$ and
$V_\gr=(\bigwedge^n V)^\sigma$.
For the Grassmann embedding of $\Gamma_n^\hform$, we have the following.
\bpr\label{prop:hermitian dual embedding}{\rm (See~\cite{De-2008a})}
The Grassmann embedding $(\emb_\gr,V_\gr)$ for $\Gamma_n^\hform$ weakly
embeds $\Gamma_n^\hform$ into $\PG(\bigwedge^n V)$; the codomain $V_\gr$
is a Baer subspace of dimension ${2n\choose n}$ over $\FF^\sigma$ that
spans $\bigwedge^n V$.
\epr

\section{Two extension results}\label{section:ext}
The main results of this section are Propositions~\ref{prop:1hyp=2hyp}~and~\ref{prop:dual to grass}. They allow us to create new
generating sets from old.
In both, we consider a non-degenerate subspace $W$ of $V$ of codimension
$2$ and Witt index $n-1$.
Proposition~\ref{prop:1hyp=2hyp} describes a relation between $\Gamma_2(W)$ and
$\Gamma_2(V)$
 and Proposition~\ref{prop:dual to grass} describes a relation between
$\Gamma_{n-1}(W)$ and $\Gamma_{n-1}(V)$.

Our first goal shall be to establish the existence of sufficiently many ``parallel
lines'' (see Lemmas~\ref{lem:Gamma2 parallel line}~and~\ref{lem:parallel lines argument}). This will be used to prove 
Corollary~\ref{cor:Gamma 2 hyperplane complement} which says that subspaces of $\Gamma_2$ have connected complements.
Proposition~\ref{prop:1hyp=2hyp} then follows easily from~Corollary~\ref{cor:Gamma 2 hyperplane complement}.
Proposition~\ref{prop:dual to grass} is independent of these results.

Let $\Theta=(\cP,\cL)$ be a partial linear space.
Let $\dist(\cdot,\cdot)$ be the distance relation on $\cP$ induced by
the natural distance in the collinearity graph of $\Theta$
 and suppose $\Theta$ has diameter $d\in \NN$.
\bde
Let $(\cD,\prec)$ be a partially ordered set  with unique minimal
element $0$.
A {\em $(\cD,\prec)$-valued distance} of $\Theta$ is a pair $(\delta,\mu)$ of maps such that the following diagram is commutative 

$$
\xymatrix{
\cP\times\cP \ar[r]^\delta \ar[dr]^{d} & \cD \ar[d]^\mu\\
 & \{0,1,\ldots,d\} \\
}
$$
and we have 
\begin{itemize}
\item[($\delta$)] for any $p,q,r\in\cP$, if $q$ and $r$ are collinear, then
  $\delta(p,q)\prec \delta(p,r)$, $\delta(p,q) = \delta(p,r)$,
 or $\delta(p,q)\succ \delta(p,r)$,
\AI{$\mu$} $\mu$ weakly preserves order: for any $\alpha,\beta\in\cD$,
$\alpha\prec \beta$ implies $\mu(\alpha)\le \mu(\beta)$,
\end{itemize}
\ede
For each $i\in \{0,1,\ldots,d\}$ and $\beta\in \cD$, define
$$\begin{array}{ll}
D_i&=\{(p,q)\mid \dist(p,q)=i\};\\
\Delta_\beta&=\{(p,q)\mid \delta(p,q)=\beta\}\\
\end{array}$$
The following important observation is easily verified.
\ble\label{lem:delta refines dist}
 The partition $\{\Delta_\beta\}_{\beta\in \cD}$ of $\cP\times\cP$
refines  $\{D_i\}_{i\in \{0,1,\ldots,d\}}$.
More precisely, $\Delta_\beta\sbe D_{\mu(\beta)}$ for all $\beta\in\cD$. 
\ele
\bde\label{dfn:gated}
We call $\Theta$ {\em gated} with respect to $(\delta,\mu)$ if, given any
point $p$ and line $l$, the set $\delta(p,l)=\{\delta(p,q)\mid q\in
l\}\sbe \cD$ has a unique minimal element.
Moreover, if  $\delta(p,l)$ has at least two elements, then the unique minimal element
 equals $\delta(p,q)$ for a unique point $q\in l$.
We then call $q$ the {\em projection} of $p$ onto $l$ and write
 $q=\proj_l(p)$.
Note here that in view of ($\delta$) we then have $\delta(p,q)\prec
\delta(p,r)$ for all $r\in l-\{q\}$.
\ede

We shall henceforth only concern ourselves with distances that are {\em
symmetric} in the sense that for any $p,q\in \cP$ we have
 $\delta(p,q)=\delta(q,p)$.
\bde\label{dfn:parallel}
Let $\Theta$ and $(\cD,\prec)$ be as above and let $(\delta,\mu)$ be a symmetric distance.
We say that two lines $l$ and $m$ in $\Theta$ are {\dfn parallel} if
\begin{itemize}
\item[(par-i)] $\delta(l,m)=\{\delta(p,q)\mid p\in l,q\in m\}$ has
exactly two elements, and
\item[(par-ii)] $\proj_m\colon l\to m$ and $\proj_l\colon m\to l$
 are mutually inverse bijections.
\end{itemize}
\ede
We now turn to the collinearity graph $\De$ of $\Gamma_2$
 and let $\dist(\cdot,\cdot)$ denote the corresponding numerical distance between points.
Since any two objects from $\Gamma$ belong to a common apartment, we can
consider
 a single apartment to determine that two points $x,y\in \Gamma_2$
 must be in one of the following two-point relations to one another:

\begin{table}[h]
\begin{centering}
\begin{tabular}{clc}
 $\delta(x,y)$ & description & $\dist(x,y)$\\
 \hline
0 &  $x=y$  & 0\\
1 &  $\langle x,y\rangle_V$ is a t.i.\ $3$-space  & 1\\
2p & $\langle x, y\rangle_V$ is a t.i.\ $4$-space & 2\\
2q & $\langle x, y\rangle_V$ is a non- t.i.\ $3$-space & 2 \\
2s & $\langle x, y\rangle_V$ is a $4$-space with a radical of dimension $2$ & 2 \\
3 &  $\langle x,y\rangle_V$ is a non-degenerate $4$-space & 3\\
\end{tabular}

\caption{}\label{table:distances}
\end{centering}
\end{table}
Here, $2{\rm p}$ stands for ``parabolic'', $2{\rm q}$ stands for ``quadrangular'', and $2{\rm s}$ stands for
``special''.
Note that relation $2{\rm p}$ only occurs if $n\ge 4$.
We now define a $(\cD,\prec)$-valued distance $(\delta,\mu)$ on $\Gamma_2$. Let
$\cD=\{0,1,2{\rm p},2{\rm q},2{\rm s},3\}$.
\begin{figure}
$$
\xymatrix{
    &    & 2{\rm p} \ar@{-}[dr] &              & \\
 0 \ar@{-}[r]& 1\ar@{-}[ur]\ar@{-}[dr] &              & 2{\rm s} \ar@{-}[r] & 3\\
    &    & 2{\rm q} \ar@{-}[ur] &              & \\
}
$$
\caption{The Hasse diagram for $(\cD,\prec)$ (rotated $90^\circ$ clockwise).}\label{fig:hasse}
\end{figure}
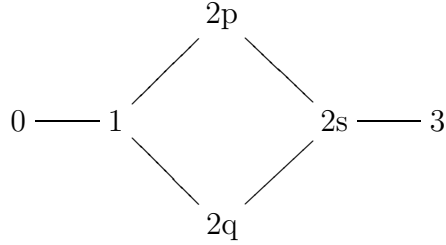
The partial order on $\cD$ is given by the Hasse diagram~in Figure~\ref{fig:hasse}.
For any two points $x,y$ we let $\delta(x,y)$ be given as in Table~\ref{table:distances} and we let $\mu$ assign $\dist(x,y)$ to $\delta(x,y)$ as given there.
It is apparent from Table~\ref{table:distances} that $(\delta,\mu)$ satisfies condition ($\mu$). Condition ($\delta$) is verified by Lemma~\ref{lem:Gatedness}.
Moreover, it is clear from middle column in Table~\ref{table:distances} that $(\delta,\mu)$ is symmetric.

\noindent
We now show that the $(\cD,\prec)$-valued distance $(\delta,\mu)$ is gated in the sense defined above.

\ble\label{lem:Gatedness}  Let $l$ be a line of $\Gamma_2$ and $p$ a
point not incident with $l.$  Then the number of relations represented
by $\{(p, q)| q \in l\}$ is at most two.  Moreover, if two
relations occur then there is a unique point $q$ on $l$ such that $(p,q)
\prec (p,r)$ for all $r \in l \setminus \{q\}.$ \ele

\bpf \nid 
Let $l$ be a line of $\Gamma_2$ and $p$ a point of $\Gamma_2.$ 
We determine the relations between $p$ and the points of $l$.
In doing
so we will identify $l$ with a flag $(P, \Pi)$ where $P$ is an isotropic
point, $\Pi$ is a totally isotropic plane and $P \subset \Pi$ and
identify $p$ with a totally isotropic line of the space $V.$  
The possible relations between $(P, \Pi)$ and $p$ are as follows:

\medskip

\nid i. $P \subset p \subset \Pi.$  In this case, $\delta(p,p)=0$ and $\delta(p,t)=1$ for all $t \in l\setminus \{p\}.$ 

\medskip

\nid ii. $P\not\subset p \subset \Pi.$  In this case, $\delta(p,t)=1$ for all $t \in l\setminus \{p\}.$ 

\medskip

\nid iii. $\Pi \cap p = P, \Pi \subset p^\perp.$  Then $\delta(p,t)=1$ for all $t \in l.$

\medskip

\nid iv. $\Pi \cap p = P, \Pi \cap p^\perp = r.$  Then $\delta(p,r) =
1$ and $\delta(p,t)=2{\rm q}$ for all $t \in l \setminus \{r\}.$

\medskip

\nid v. $\Pi \cap p = Q \ne P, \Pi \subset p^\perp.$  If $r = \langle P ,Q\rangle_V$
then $\delta(p,r) = 1$ and $\delta(p,t)=2{\rm p}$ for all $t \in l
\setminus \{r\}.$

\medskip

\nid vi. $\Pi \cap p = Q \ne P, \Pi \cap p^\perp = \langle P, Q \rangle_V= r.$  Then
$\delta(p,r) = 1$ and $\delta(p,t)=2{\rm s}$ for all $t \in l
\setminus \{r\}.$

\medskip

\nid vii. $\Pi \cap p = Q \ne P, P \nsubseteq \Pi \cap p^\perp.$  If $r =
\langle P,Q\rangle_V$ then $\delta(p,r)=2{\rm q}$ and $\delta(p,t)=2{\rm s}$
for all $t \in l \setminus \{r\}.$

\medskip

\nid viii. $\Pi \cap p = 0, \Pi \subset p^\perp.$  Then $\delta(p,t)=2{\rm p}$ for all $t \in l.$

\medskip

\nid ix. $\Pi \cap p = 0, \Pi \cap p^\perp = r \in l.$  Then $\delta(p,r)=2{\rm p}$ and $\delta(p,t)=2{\rm s}$ for all $t \in l
\setminus \{r\}.$

\medskip

\nid x. $\Pi \cap p = 0, \Pi \cap p^\perp = x$ is a point of $\Gamma_2$
not on $l.$   Now $\delta(p,t)=2{\rm p}$ for all $t \in l.$

\medskip

\nid xi. $\Pi \cap p = 0, \Pi \cap p^\perp = P.$  Now $\delta(p,t)=2{\rm s}$ for all $t \in l.$

\medskip

\nid xii. $\Pi \cap p = 0, \Pi \cap p^\perp = Q \ne P.$ If $r = \langle P, Q\rangle_V$
then $\delta(p,r)=2{\rm s}$ and $\delta(p,t)=3$ for all 
$t \in l \setminus \{r\}.$
\epf

Since by Lemma~\ref{lem:Gatedness}, $\Gamma_2^\form$ is gated in the sense of Definition~\ref{dfn:gated}, it makes sense to talk about parallel lines.

\ble\label{lem:Gamma2 parallel line}
Given any two non-collinear points $p$ and $q$, there exist
parallel lines $l$ and $m$ on $p$ and $q$ so that $\proj_m(p)\ne q$.
\ele
\bpf
Let $\cE=\{e_i,f_i\mid i=1,2,\ldots,n\}$ be a hyperbolic basis
 such that $p,q\in \Sigma(\cE)$ with $p=\langle e_1,e_2\rangle$.
Now $l=(p\cap p',\langle p,p'\rangle)$
 and $m=(q\cap q',\langle q,q'\rangle)$ are given by the following table, where $q'=\proj_m(p)$ and $p'=\proj_l(q)$.
$$\begin{array}{c|c|c|c|c}
 q  & p' & q' & \delta(p,q)= \delta(p',q') & \delta(p,q')=\delta(p',q)  \\
\hline
\langle f_1,f_2\rangle & \langle e_1,e_3\rangle & \langle f_1,f_3\rangle & 3 & 2{\rm s}\\
\hline
\langle f_1,f_3\rangle & \langle e_1,f_3\rangle  & \langle f_1,e_2\rangle & 2{\rm s} & 2{\rm q}\\
\hline
\langle f_1,e_2\rangle & \langle e_2,f_3\rangle  & \langle e_2,e_3\rangle  & 2{\rm q} & 1\\
\hline
\langle e_3,e_4\rangle & \langle e_1,e_3\rangle & \langle e_2,e_4\rangle & 2{\rm p} & 1\\
\end{array}$$
%
%
\epf
 
Note that in Lemma~\ref{lem:Gamma2 parallel line} we do not require
that, given any two points $p$ and $q$ {\em and a line $l$ on $p$},
there exists a line $m$ on $q$ that is parallel to $l$. This is not true
even in this geometry. However, we do not need that for proving Lemma~\ref{lem:parallel lines
argument}.

\medskip
\noindent
The first author learned the following useful argument from Andries Brouwer.

\ble\label{lem:parallel lines argument}
Let $\Theta$ be a thick partial linear space with point-set $\cP$ of
finite diameter whose two-point relations
 can be ordered in some way so that given any two points $p$ and $q$,
there exist parallel lines $l$ and $m$ on $p$ and $q$ so that
$\proj_m(p)\ne q$.
Then the complement of any proper subspace is connected.
\ele
\bpf
Let $\cH$ be a proper subspace and set $\cA=\cP-\cH$.

For any two points $p,q\in \cA$ we show that there is a path of
 points and lines in $\cA$ from $p$ to $q$.

If $p$ and $q$ are on a common line $l$, then $p,l,q$ is a path
 in $\cA$ connecting $p$ and $q$ and we are done.

Now let $p$ and $q$ be non-collinear, and let $l$ and $m$ be the
 parallel lines assumed to exist.
Then because lines of $\Theta$ are thick and intersect the
subspace $\cH$ in none, one, or all of their points,
 one of the following must happen:
 (1) $r=\proj_m(p)\in \cA$,
 (2)  $s=\proj_l(q)\in \cA$,
 or
 (3) there exist points $t, u\in \cA$ different from
 $p$, $q$, $r$, $s$ and on $l$ and $m$, respectively, such that
  $\proj_m(t)=u$.

In all cases, $p$ and $q$ are collinear to points of $\cA$ that are
nearer to each other than are $p$ and $q$ (since $\proj_m(p)\ne q$).
Since the diameter of $\Theta$ is finite, by repeating this argument we
find a path in $\cA$  connecting $p$ to $q$.
\epf

\medskip

From Lemmas~\ref{lem:Gamma2 parallel line}~and~\ref{lem:parallel lines argument}, we obtain the following result.

\bco\label{cor:Gamma 2 hyperplane complement}
The complement of a proper subspace of $\Gamma_2$ is connected.
\eco

\bpr\label{prop:1hyp=2hyp}
Assume $n\ge 3$ and let $W$ be a non-degenerate $2(n-1)$-space  of Witt
index $n-1$ in $V$. Then
 $\cH=\{K\in\cP({\Gamma_2^\form})\mid K\cap W\ne \{0\}\}$
is a hyperplane of $\Gamma^\form_2$.
This is a maximal subspace.
\epr

\bpf
Let the line $l$ of $\Gamma_2^\form$ consist of all $2$-objects $K$ with
 $A\sbe K\sbe C$ for some $1$-object $A$ and $3$-object $C$.
Now $\langle A,A^\perp\cap W\rangle/A$ is a hyperplane or all of the
polar space $A^\perp/A$.
It follows that, accordingly, either one or all of the points of $l$
meet $W$, which proves the lemma.
The latter claim follows from Corollary~\ref{cor:Gamma 2 hyperplane
complement}.
\epf

\medskip\noindent
Proposition~\ref{prop:dual to grass} can be used for an inductive argument on the
 Witt index $n$ of $V$.
\bpr\label{prop:dual to grass}
Assume $\dim V\ge 6$.
Let $W\le V$ be a non-degenerate subspace of codimension $2$ and Witt
index $n-1$, and let $P,Q$ be two isotropic points in the hyperbolic line
$W^\perp$.
If $S$ is a set of points in $\Gamma_{k}(W)$, then
 $$\langle S\rangle_{\Gamma_{k}(W)}\sbe \langle S,
\Res_{k}(P),\Res_{k}(Q)\rangle_{\Gamma_{k}(V)}.$$
\epr
\bpf
Clearly points of $\Gamma_k(W)$ are also points of $\Gamma_k(V)$.
If $k<n-1$, then also lines of the former are lines of the latter.
In that case the result is trivial.

Now let $k=n-1$. Note that $W$ has Witt index $(n-1)$ and $V$ has Witt
index $n$.
Let $x,y\in \Gamma_{n-1}(W)$ be on some line $l$ of $\Gamma_{n-1}(W)$.
Then there is a t.i.\ $(n-2)$-space $L\le W$ such that $l$ consists of
all t.i.\
$(n-1)$-spaces in  $L^\perp\cap W$.
Note that $L^\perp\cap W$ is a non t.i.\ $n$-space and so the points on
$l$ are not collinear in $\Gamma_{n-1}(V)$.

However, all points of $l$ are contained in the residue
$\Res_{\Gamma_{n-1}}(L)\cong \Gamma_1(L^\perp/L)$.
Thus, by  Proposition~\ref{prop:Gamma_1 spanned by apartment} this residue is
generated as a subspace of $\Gamma_{n-1}(V)$ by
 the set $\{x,y,\langle L,P\rangle_V,\langle L,Q\rangle_V\}$.
We are done since $\langle L,P\rangle\in\Res_{n-1}(P)$ and
 $\langle L,Q\rangle\in\Res_{n-1}(Q)$.
 \epf
\section{Generating sets}\label{section:gensets}
We inductively define the following set
$\cS^\form_{n,k}=\cS^\form_{n,k}(V)$ of $k$-objects in $\Gamma^\form_k$:
Consider a non-degenerate $2(n-1)$-space $W$ of Witt index $n-1$ of $V$
and let $P$ and $Q$ form a hyperbolic pair
 in $W^\perp$.
We formally define $\cS_{n,0}(V)=\{0\}$, thus containing one element.
Let $\cS_{n,1}(V)$ and $\cS_{n,n}(V)$ be the generating sets
 of $\Gamma_1^\form$ and $\Gamma_n^\form$ respectively.
The existence of these sets is guaranteed by Proposition~\ref{prop:Gamma_1 spanned by apartment} and Theorem~\ref{thm:bruce's theorem}.
Now assume $2\le k\le n-1$.
Using induction on $n$, we choose $\cS^\form_{n-1,l}(W)$ for
$l=k,k-1,k-2$.
Let $$\hcS^\form_{n-1,k-1}(W)=\{\langle L,P\rangle, \langle
L,Q\rangle\mid L\in \cS_{n-1,k-1}^\form(W)\}.$$

Finally we describe the set $\hcS^\form_{n-1,k-2}(W)$.
Let $M$ be any $(k-2)$-object in $W$. Then, the collection of all
$k$-objects $K$ containing $M$ and intersecting $M^\perp\cap W$
non-trivially, is a hyperplane in the residue of $M$ in the $k$-Grassmannian
 (compare Proposition~\ref{prop:1hyp=2hyp}).
Let $\what{M}$ be some $k$-object on $M$, not contained in this hyperplane.
Now set
$$
\hcS^\form_{n-1,k-2}(W)=\{\what{M}\mid M\in\cS_{n-1,k-2}^\form(W)\}.
$$
Using induction on $n$ we now define
$$
\cS^\form_{n,k}(V)=\cS_{n-1,k}^\form(W)\uplus
\hcS_{n-1,k-1}^\form(W)\uplus \hcS_{n-1,k-2}^\form(W).
$$
We note that in case $k=2$, $\hcS^\form_{n-1,k-2}$ contains one element.

\ble\label{lem:binomial}
\quad\begin{itemize}
\AI{a} For $1\le k\le m$ integers, we have
${m\choose k}={m-2\choose k}+2{m-2\choose k-1}+{m-2\choose k-2}$.
\AI{b} For $1\le k\le n$ integers, we have
$$\begin{array}{ll}
{2n\choose k}&={2(n-1)\choose k}+2{2(n-1)\choose k-1}+{2(n-1)\choose
k-2},\mbox{ and } \\
 {2n\choose k}-{2n\choose k-2}&= {2(n-1)\choose k}+2{2(n-1)\choose
k-1}-2{2(n-1)\choose k-3}-{2(n-1)\choose k-4}.\\
\end{array}$$
\end{itemize}
\ele
\bpf
(a) Apply the binomial theorem to  $(x+1)^m=(x+1)^2\cdot (x+1)^{m-2}$ and
compute the coefficient of $x^k$.
(b) This is immediate from (a).
\epf

\bco\label{cor:size of gensets}
Let $1\le k\le n$ be integers. Then,
\begin{itemize}
\AI{a}
$|\cS^\sform_{n,k}|={2n\choose k}-{2n\choose k-2}$,
\AI{b}
$|\cS^\hform_{n,k}|={2n\choose k}$.
\end{itemize}
\eco
\bpf
For $k=1,n$, this is clear by 
Proposition~\ref{prop:Gamma_1 spanned by
apartment} and Theorem~\ref{thm:bruce's theorem}.
Now let $k\ge 2$.
By construction, we have
$|\cS^\form_{n,k}|=|\cS^\form_{n-1,k}|+2|\cS^\form_{n-1,k-1}|+|\cS^\form_{n-1,k-2}|$.
Therefore the lemma follows by induction on $2\le k \le n-1$ and for
each $k$ by induction on $n>k$ from Lemma~\ref{lem:binomial}.
\epf

\medskip
We recall the following definition from~\cite{BlPa2001}.
Given a set $\cS$ of points in $\Gamma_k^\form$, we call an object $X$
of $\Gamma$ {\dfn $\cS$-full}
 if $\Res_k(X)\le \langle \cS\rangle_{\Gamma_k^\form}$.

\bpr\label{prop:span of gensets}
Let $1\le k\le n$ be integers. Then,
$$\langle \cS^\form_{n,k}\rangle_{\Gamma^\form_k}={\Gamma^\form_k}.$$
\epr
\bpf
The case $k=1$ follows from Proposition~\ref{prop:Gamma_1 spanned by apartment}.
The case $k=n$ follows from Theorem~\ref{thm:bruce's theorem}.
We now continue by induction on $n$.
Set $\cS=\cS^\form_{n,k}(V)$.
Recall $W$ is a non-degenerate $2(n-1)$-dimensional subspace of $V$ of
Witt index $n-1$ and let $P$ and $Q$ be two points of the hyperbolic
line $W^\perp$.
Let $\cS_l(W)=\cS^\form_{n-1,l}(W)$.

We first note that the points $P$ and $Q$ are $\cS$-full.
Let us see why this is so.
For each $L\in\cS_{n-1,k-1}(W)$, we have $\langle L, P\rangle\in\cS$.
Now let $L_1,L_2$ be $(k-1)$-objects that are on some line of
$\Gamma_{k-1}^\form(W)$.
Note that, since $k-1<n-1$, this line is of the form $(A,B)$ for some
$(k-2)$-object $A$ and
 $k$-object $B$ with $A\le B\le W$.
Then, $(\langle A,P\rangle, \langle B,P\rangle)$ is a line of
${\Gamma^\form_k(V)}$ containing
 the points $\langle L_1,P\rangle$ and $\langle L_2,P\rangle$.
It follows that for each $L\in\Gamma_{k-1}^\form(W)=\langle
\cS_{k-1}(W)\rangle$ we have
 $\langle L,P\rangle\in\langle \cS\rangle_{\Gamma^\form_k(V)}$.
In particular, $P$ and, similarly, $Q$ are $\cS$-full.
It now follows from Proposition~\ref{prop:dual to grass} that
  $$\langle \cS_k(W)\rangle_{\Gamma_k^\form(W)}\sbe \langle
\cS_k(W)\rangle_{\Gamma^\form_k(V)}.$$
We now show that every $(k-1)$-object $L\in\Gamma_{k-1}^\form(W)$ is
$\cS$-full.
The residue of $L$ is of type $\Gamma^\form_1(L^\perp/L)$.
Now by the preceding, $\langle \cS\rangle_{\Gamma_k^\form(V)}$ contains all points of
$\Gamma^\form_1(L^\perp\cap W/L)$
as well as the points $\langle P,L\rangle$ and $\langle Q,L\rangle$. Therefore, by the case $k=1$
(Proposition~\ref{prop:Gamma_1 spanned by apartment})
 $\langle \cS\rangle$ contains all points of $\Gamma^\form_1(L^\perp/L)$.

We'll now show that every $(k-2)$-object $M\in\Gamma_{k-2}^\form(W)$ is
$\cS$-full.
The residue of $M$ is of type $\Gamma^\form_2(M^\perp/M)$.
Now by the preceding, $\langle \cS\rangle_{\Gamma_k^\form(V)}$ contains all $2$-objects of
$\Gamma^\form_2(M^\perp/M)$
 meeting $(M^\perp\cap W)/M$ non-trivially.
Therefore, by Proposition~\ref{prop:1hyp=2hyp}, it suffices to show that
 $\langle \cS\rangle_{\Gamma_k^\form(V)}$ also contains one $k$-object $K$ with $M=K\cap W$.
For each $M\in \cS_{k-2}(W)$, this is true by definition of $\cS$.
Note that this settles the $k=2$ case.

 From now on we may assume $k\ge 3$.
Suppose that $M\in \Gamma_k^\form(W)\setminus\cS_{k-2}(W)$.
Now let $M_1,M_2$ be $\cS$-full $(k-2)$-objects that are on some line of
$\Gamma_{k-2}^\form(W)$ containing $M$.
Note that, since $k-2<n-1$, this line is of the form $(A,B)$ for some
$(k-3)$-object $A$ and
 $(k-1)$-object $B$ with $A\le B\le W$.
Since $k\le n-1$, there exists a $(k+1)$-object of the form
 $\langle B,R\rangle$ where $R$ is of dimension $2$ and disjoint from $W$.
Now $(\langle A,R\rangle, \langle B,R\rangle)$ is a line of
${\Gamma^\form_k(V)}$ spanned by the points
 $\langle M_1,R\rangle$ and $\langle M_2,R\rangle$.
Since these points belong to $\langle \cS\rangle_{\Gamma_k^\form(V)}$, so does $\langle
M,R\rangle$.
Since $M=\langle M, R\rangle\cap W$, we find that $M$ is $\cS$-full.

Since, for each $l=k-2,k-1,k$, every $l$-object in $W$ is $\cS$-full,
and every $k$-object intersects $W$
 in an object of type $k-2$, $k-1$, or $k$, we find that
${\Gamma^\form_k(V)}=\langle \cS\rangle_{\Gamma_k^\form(V)}$.
\epf

\bpr\label{prop:genset independent}
The image of $\cS^\form_{n,k}(V)$ under $e_\gr$ forms a basis for $V_\gr$.
\epr
\bpf
Let $I=\{1,2,\ldots,2n\}$.
Given any basis $\{e_i\mid i\in I\}$ for $V$, it is known that
$\bigwedge^k V$ has basis
 $\{e_K\mid K\sbe I,\ |K|=k\}$, where
  $e_K=\wedge_{k\in K} e_k$ is taken in order of increasing $k$.
Let $W$, $P$ and $Q$ be as in the construction of $\cS_{n,k}^\form(V)$.
Choosing our basis so that $\{e_i\mid i\in I-\{2n-1,2n\}\}$ is a basis
for $W$ and
 $P=\langle e_{2n-1}\rangle$ and $Q=\langle e_{2n}\rangle$ it follows that
\beq\label{eqn:wedgek decomp}
\bigwedge^k V=\bigwedge^k W\oplus \left(\bigwedge^{k-1} W\wedge
P\right)\oplus\left(\bigwedge^{k-1} W\wedge Q\right)\oplus
\left(\bigwedge^{k-2} W\wedge P\wedge Q\right).
\eeq
We wish to show that for all $n\ge 2$ and $1\le k\le n$, the images of
$\cS_{n,k}^\form$ under the embedding $\emb_\gr$ are linearly independent.
We note that in these cases the geometry $\Gamma_k^\form$ is embedded
into a subspace of $\bigwedge^k V$.
Therefore it suffices to show that these images are linearly independent
in $\bigwedge^k V$.
With slight abuse of language, we shall say that a set $\cS$ of points
of $\Gamma_k^\form$ is ``$\emb_\gr$-independent''
 if the set $\emb_\gr(\cS)$ is linearly independent in $\bigwedge^k V$.

We show that $\cS_{n,k}^\form(V)$ is $\emb_\gr$-independent, using
induction on $n$ and, for each $n$ we distinguish cases
 $k=1$, $k=n$ and $2\le k\le n-1$.
If $k=1$ and $n\ge 1$, then $V_\gr=V$ and $\cS_{n,1}^\form$ is simply
the set of
 $1$-spaces spanned by a hyperbolic basis for $V$, so the claim holds.

Next, we address the special case $n=k$.
If $\form=\sform$, we note that $\Gamma_n^\sform$ is naturally embedded
into $\bigwedge^n V$ and
 so  the claim follows directly from Theorem~\ref{thm:bruce's theorem}.
 From Proposition~\ref{prop:hermitian dual embedding} we know that the
geometry $\Gamma_n^\hform$ is weakly, but not fully embedded into
$\bigwedge^n V$ since
 this geometry and its embedding are naturally defined over the fixed
field $\FF^\sigma$ of $\FF$ under $\sigma$.
However, we can still view the image of $\cS_{n,n}^\hform$ under $e_\gr$
as a set of vectors of $\bigwedge^n V$.
We then note that
$V_\gr\otimes_{\FF^\sigma}\FF=\bigwedge^n V$,
 so in particular, by Theorem~\ref{thm:bruce's theorem}, the image
$\cS_{n,n}^\hform$ under $e_\gr$ spans $\bigwedge^n V$.

If we view $\cS_{n-1,n-1}^\hform(W)$ as a subset of $\cS_{n,n-1}(V)$,
then its image under
 $\emb_\gr$ is linearly independent in $V_\gr$ by
Theorem~\ref{thm:bruce's theorem}.
By the preceding paragraph, $\emb_\gr(\cS_{n-1,n-1}^\hform(W))=\{\langle
\emb_\gr(s)\rangle_\FF\mid s\in \cS_{n-1,n-1}^\hform(W)\}$.
And since the codomain $W_\gr$ of the embedding $e_\gr$ of
$\Gamma_{n-1}(W)$ linearly spans $\bigwedge^{n-1} W$, the spanning set
$\cS_{n-1,n-1}^\hform(W)$ is $\emb_\gr$ independent.
This settles the base case $k=n$.

Next, assume $2\le k\le n-1$. We note that
$$\cS^\form_{n,k}(V)=\cS^\form_{n-1,k}(W)\uplus\hcS^\form_{n-1,k-1}(W)\uplus\hcS^\form_{n-1,k-2}(W).$$
Now assume that $\cS_{n-1,k}^\form(W)$, $\cS_{n-1,k-1}^\form(W)$ and
 $\cS_{n-1,k-2}^\form(W)$ are all independent.
Note that
$$
\begin{array}{ll}
\emb_\gr(\cS_{n-1,k}^\form(W))&\sbe \bigwedge^k W,\\
\emb_\gr(\hcS_{n-1,k-1}^\form(W))&\sbe \left( \bigwedge^{k-1}W\wedge
P\right)\oplus
\left(\bigwedge^{k-1}W \wedge Q\right),\\
\emb_\gr(\hcS_{n-1,k-2}^\form(W))&\sbe \bigwedge^{k-2}W\wedge P\wedge Q.\\
\end{array}
$$
Since $\cS^\form_{n-1,k-1}(W)$ is independent, so is the image of
$\hcS_{n-1,k-1}^\form(W)$ in $\bigwedge^{k-1}W\wedge P$, and likewise
for $Q$.
For the same reason, the image of $\hcS_{n-1,k-2}^\form(W)$ in
$\bigwedge^{k-2} W\wedge P\wedge Q$ is independent.
The direct sum decomposition~(\ref{eqn:wedgek decomp}) shows that the
set $\hcS^\form_{n-1,k-1}(W)$ is independent, and that the sets
$\cS_{n-1,k}^\form(W)$, $\hcS_{n-1,k-1}^\form(W)$, and 
$\hcS_{n-1,k-2}^\form(W)$ are also pairwise independent.
\epf

\medskip
%
\bpr\label{prop:symplectic grassmann embedding}
Let $1\le k\le n$.
The Grassmann embedding $(\emb_\gr,V_\gr)$ for $\Gamma_k^\sform$ has
dimension ${2n\choose k}-{2n\choose k-2}$.
\epr
\bpf
In Proposition~\ref{prop:genset independent} we showed directly that
the image of the  generating set $\cS^\sform_{n,k}(V)$ forms a basis for
$V_\gr$.
It then follows from Corollary~\ref{cor:size of gensets} that
$\dim(V_\gr)=|\cS^\sform_{n,k}(V)|={2n\choose k}-{2n\choose k-2}$.
\epf


\begin{thebibliography}{10}
\expandafter\ifx\csname url\endcsname\relax
  \def\url#1{\texttt{#1}}\fi
\expandafter\ifx\csname urlprefix\endcsname\relax\def\urlprefix{URL }\fi

\bibitem{BeSh2004}
C.~D. Bennett and S.~Shpectorov.
\newblock A new proof of a theorem of {P}han.
\newblock {\em J. Group Theory}, 7(3):287--310, 2004.


\bibitem{Bl1999}
R.~J. Blok, On geometries related to buildings, Ph.D. thesis, Delft
University
  of Technology, supervisor: Prof. Dr. A.E. Brouwer (1999).

\bibitem{Bl2003}
R.~J. Blok, The generating rank of the symplectic line-{G}rassmannian,
  Beitr\"age Algebra Geom. 44~(2) (2003) 575--580.

\bibitem{Bl2007}
R.~J. Blok, The generating rank of the symplectic grassmannians:
Hyperbolic and
  isotropic geometry, Europ. J. Comb. 28~(5) (2007) 1368--1394.

\bibitem{Blo2010}
R.~J. Blok.
\newblock Highest weight modules and polarized embeddings of shadow spaces.
\newblock 2010.


\bibitem{BlBr1998a}
R.~J. Blok, A.~E. Brouwer, \protect{Spanning point-line geometries in
buildings
  of spherical type}, J. Geom. 62 (1998) 26--35.

\bibitem{BlPa2001}
R.~J. Blok, A.~Pasini, Point-line geometries with a generating set that
depends
  on the underlying field, in: Finite geometries, Vol.~3 of Dev. Math.,
Kluwer
  Acad. Publ., Dordrecht, 2001, pp. 1--25.

\bibitem{BrCoNe1989}
A.~Brouwer, A.~Cohen, A.~Neumaier, Distance-regular graphs, no.~18 in
  Ergebnisse der Mathematik und ihrer Grenzgebiete 3, Springer, Berlin,
1989.

\bibitem{BuSh1974}
F.~Buekenhout and E.~Shult.
\newblock On the foundations of polar geometry.
\newblock {\em Geometriae Ded.}, 3:155--170, 1974.


\bibitem{CohCoo1998}
A.~M. Cohen and B.~N. Cooperstein.
\newblock Lie incidence systems from projective varieties.
\newblock {\em Proc. Amer. Math. Soc.}, 126(7):2095--2102, 1998.



\bibitem{Co1997}
B.~N. Cooperstein, On the generation of dual polar spaces of unitary
type over
  finite fields, European J. Combin. 18~(8) (1997) 849--856.

\bibitem{Co1998}
B.~N. Cooperstein, On the generation of dual polar spaces of symplectic type
  over finite fields, J. Combin. Theory Ser. A 83~(2) (1998) 221--232.

\bibitem{Co2003}
B.~N. Cooperstein, Generation of embeddable incidence geometries: a
survey, in:
  Topics in diagram geometry, Vol.~12 of Quad. Mat., Dept. Math.,
Seconda Univ.
  Napoli, Caserta, 2003, pp. 29--57.

\bibitem{CoSh1997}
B.~N. Cooperstein, E.~E. Shult, Frames and bases of {L}ie incidence
geometries,
  J. Geom. 60~(1-2) (1997) 17--46.

\bibitem{DavPri2002}
B.~A. Davey and H.~A. Priestley.
\newblock {\em Introduction to lattices and order}.
\newblock Cambridge University Press, New York, second edition, 2002.

\bibitem{De-2008a}
B.~De~Bruyn, On the {G}rassmann-embeddings of the {H}ermitian dual polar
  spaces, Linear Multilinear Algebra 56~(6) (2008) 665--677.

\bibitem{De-2009}
B.~De~Bruyn.
\newblock Some subspaces of the {$k$}th exterior power of a symplectic vector
  space.
\newblock {\em Linear Algebra Appl.}, 430(11-12):3095--3104, 2009.


\bibitem{DePa2007}
B.~De~Bruyn, A.~Pasini, Generating symplectic and {H}ermitian dual polar
spaces
  over arbitrary fields nonisomorphic to {$\mathbb F\sb 2$}, Electron. J.
  Combin. 14~(1) (2007) Research Paper 54, 17 pp. (electronic).


\bibitem{Di1955}
J.~Dieudonn\'{e}, La g\'{e}om\'{e}trie des groupes classiques, Vol.~5 of
  Ergebnisse der Mathematik und ihrer Grenzgebiete, Springer, Berlin, 1955.
 
 \bibitem{Fu1997}
W.~Fulton.
\newblock {\em Young Tableaux}, volume~35 of {\em London Math. Soc. Student
  Texts}.
\newblock Cambridge University Press, Cambridge, 1997.


\bibitem{KuSt2004}
H.~Kurzweil and B.~Stellmacher.
\newblock {\em The theory of finite groups}.
\newblock Universitext. Springer-Verlag, New York, 2004.
\newblock An introduction, Translated from the 1998 German original.



\bibitem{Li2002}
P.~Li.
\newblock On the universal embedding of the {$U\sb {2n}(2)$} dual polar space.
\newblock {\em J. Combin. Theory Ser. A}, 98(2):235--252, 2002.

\bibitem{Sh1992}
E.~Shult.
\newblock \protect{Geometric Hyperplanes of Embeddable Grassmannians}.
\newblock {\em J. Algebra}, 145:55--82, 1992.

\bibitem{St1997}
R.~P. Stanley.
\newblock {\em Enumerative Combinatorics, Volume I}, volume~49 of {\em
  Cambridge studies in advanced mathematics}.
\newblock Cambridge University Press, 1997.


\bibitem{Ta1992}
D.~E.~Taylor, The geometry of the classical groups, 1st Edition, Vol.~9 of
  Sigma Series in Pure Math., Heldermann, Berlin, 1992.

\bibitem{Ti1974}
J.~Tits.
\newblock {\em Buildings of spherical type and finite {BN}-pairs}.
\newblock Springer-Verlag, Berlin, 1974.
\newblock Lecture Notes in Mathematics, Vol. 386.

\bibitem{Vel1959}
F.~D. Veldkamp.
\newblock Polar geometry. {I}, {II}, {III}, {IV}, {V}.
\newblock {\em Nederl. Akad. Wetensch. Proc. Ser. A 62; 63 = Indag. Math. 21
  (1959), 512-551}, 22:207--212, 1959.

\end{thebibliography}
\end{document}